\documentclass[english]{article}
\usepackage[T1]{fontenc}
\usepackage[latin9]{inputenc}
\usepackage{geometry}
\usepackage{verbatim}
\usepackage{amsmath}
\usepackage{amssymb}
\usepackage{esint}

\makeatletter
\newcommand{\lyxaddress}[1]{
	\par {\raggedright #1
	\vspace{1.4em}
	\noindent\par}
}

\makeatother

\makeatother

\usepackage{babel}

\makeatother

\usepackage{babel}

\makeatother

\usepackage{babel}
\begin{document}

\title{The Fourier-Legendre series of Bessel functions of the first kind
and the summed series involving $\,_{2}F_{3}$ hypergeometric functions
that arise from them}

\author{Jack C. Straton}
\maketitle


\lyxaddress{Department of Physics, Portland State University, Portland, OR 97207-0751,
USA}

\lyxaddress{straton@pdx.edu}
\begin{abstract}
The Bessel function of the first kind $J_{N}\left(kx\right)$ is expanded
in a Fourier-Legendre series, as is the modified Bessel functions
of the first kind $I_{N}\left(kx\right)$. The purpose of these expansions
in Legendre polynomials was not an attempt to rival established \emph{numerical
methods} for calculating Bessel functions, but to provide a form for
$J_{N}\left(kx\right)$ useful for \emph{analytical} work in the area
of strong laser fields, where analytical integration over scattering
angles is essential. Despite their primary purpose, we can easily
truncate the series at 21 terms to provide 33-digit accuracy that
matches IEEE extended precision in some compilers. The analytical
theme is furthered by showing that infinite series of like-powered
contributors (involving $\,_{2}F_{3}$ hypergeometric functions) extracted
from the Fourier-Legendre series may be summed, having values that
are inverse powers of the eight primes $1/\left(2^{i}3^{j}5^{k}7^{l}11^{m}13^{n}17^{o}19^{p}\right)$
multiplying powers of the coefficient $k$.
\end{abstract}
\vspace{2pc}
 \textit{Keywords}: Bessel functions, Fourier-Legendre series, Laplace
series, Polynomial approximations, Computational methods \\
 \textit{2020 Mathematics Subject Classification}: 33C10, 42C10, 41A10,
33F10, 65D20, 68W30
 \\

\section{Introduction}

\noindent Bessel functions find application in countless fields, so
there has naturally been a great deal of research into how best to
calculate them. Just a few of these approaches range from polynomial
approximations,\cite{AS p. 369 No. 9.4.1,PolyMillane} to solving
integral equations,\cite{int_Gross} to rational approximations,\cite{rat_Maass,rat_Mart=0000EDn}
to expansions in Jacobi or Chebyshev polynomials, \cite{jacobi_cheb_Wimp}.
The latter are particularly prized by many workers.\cite{cheb_Clenshaw,Cheb_Zhang,cheb_Coleman}

However, the present paper has an entirely different goal from \emph{numerical}
approximations to Bessel functions. We seek a representation of the
Bessel function $J_{N}\left(kx\right)$ in a form that will be useful
for \emph{analytical} integration over angular variables. The work
most closely related to the present approach are trigonometric approximations.\cite{trig_Molinari}

This need arises from the use of the Strong Field Approximation (SFA)
\cite{Reiss1,Reiss3,Reiss4,Reiss5,Reiss6,Faisel2} to calculate atomic
transition amplitudes in intense laser fields. Unlike perturbation
expansions that will not converge if an applied laser field is large,
the SFA is an analytical approximation that is non-perturbative. Keating
\cite{KeatingPhD} applied it to the production of the positive antihydrogen
ion and found that the reduction of the transition amplitudes to analytic
form required the expansion of the resultant Bessel functions $J_{N}\left(kx\right)$
in a series of spherical harmonics. We have been unable to discover
such a Laplace series \cite{Haber} in the literature, nor does there
seem to be an expansion of the Bessel function in a series of Legendre
polynomials, to which the Laplace series reduces in the common case
where the function is independent of the azimuthal angle. We recreate
herein Keating's derivation, stripped of the specialized SFA terminology,
and find that we can recast his coefficients from another layer of infinite series to a $\,_{2}F_{3}$ generalized
hypergeometric function. We also extend the method to the modified
Bessel function of the first kind $I_{N}\left(kx\right)$.

Though designed as an analytical tool, this Fourier-Legendre series
is easily converted to a series in powers, we do so as a check on
its validity. Computer code is given in the appendix for these series that provide 33-digit accuracy, matching IEEE extended precision in some compilers.

Of further analytical significance,  like-powered contributors extracted
from the Fourier-Legendre series may be summed, having values that
are inverse inverse powers of the eight primes $1/\left(2^{i}3^{j}5^{k}7^{l}11^{m}13^{n}17^{o}19^{p}\right)$
multiplying even powers of the coefficient $k$.

\section{The Fourier-Legendre series of a Bessel function of the first kind}

We begin with the assumption \cite{Kellogg}that the series 
\begin{equation}
J_{N}\left(kx\right)=\sum_{L=0}^{\infty}a_{LN}\left(k\right)P_{L}(x)\quad\label{eq:Fourier-Legendre series}
\end{equation}
 converges uniformly,%
\footnote{Let D be a region in which the above series converges for each value
of x. Then the series can be said to converge uniformly in D if, for
every $\varepsilon>0$, there exists a number $N'(\varepsilon)$
such that, for $n>N'$, it follows that

\[
\left. \left|J_{N}\left(kx\right)-\sum_{L=0}^{n}a_{LN}\left(k\right)P_{L}(x)\right|= \left|\sum_{L=n+1}^{\infty}a_{LN}\left(k\right)P_{L}(x)\right| \right. <\varepsilon
\]

for all x in D.
} where the coefficients are given by the orthogonality of the Legendre
polynomials,

\begin{equation}
a_{LN}\left(k\right)=\frac{2L+1}{2}\int_{-1}^{1}J_{N}\left(kx\right)P_{L}(x)dx\;.\label{eq:coefficients}
\end{equation}

To find these coefficients we following Keating's approach, first
using Heine's integral representation of the Bessel function \cite{MagnusOberhettingerSoni}
for integer indices

\begin{equation}
J_{N}(kx)=\frac{i^{-N}}{\pi}\int_{0}^{\pi}e^{ikx\cos\theta}\cos\left(N\theta\right)d\theta\;,\label{eq:jint}
\end{equation}
 so that

\begin{equation}
a_{LN}\left(k\right)=\frac{2L+1}{2}\int_{-1}^{1}\left[\frac{i^{-N}}{\pi}\int_{0}^{\pi}e^{ikx\cos\theta}\cos\left(N\theta\right)d\theta\right]P_{L}(x)dx\;.\label{eq:curleyjxtheta}
\end{equation}
 By switching the order of integration,

\begin{equation}
a_{LN}\left(k\right)=\frac{2L+1}{2}\frac{i^{-N}}{\pi}\int_{0}^{\pi}\left[\int_{-1}^{1}e^{ikx\cos\theta}P_{L}(x)dx\right]\cos\left(N\theta\right)d\theta\;,\label{eq:scriptjbeforelpinsertion}
\end{equation}
 we can expand the plane wave in a series of products of Spherical
Bessel functions $j_{m}(y)$ and Legendre polynomials\cite{GR5 p. 987 No. 8.511.4,Joachain p. 671 Eq. (B.46)}

\begin{equation}
e^{ikx\cos\theta}=\sum_{l'}(2l'+1)i^{l'}j_{l'}\left(k\cos\theta\right)P_{l'}(x)\;.\label{eq:GR5 p. 987 No. 8.511.4}
\end{equation}

\noindent Then

\begin{equation}
\begin{array}{ccc}
a_{LN}\left(k\right) & = & \frac{2L+1}{2}\frac{i^{-N}}{\pi}\int_{0}^{\pi}\left[\int_{-1}^{1}\left(\sum_{l'=0}^{\infty}(2l'+1)i^{l'}j_{l'}\left(k\cos\theta\right)P_{l'}(x)\right)P_{L}(x)dx\right]\cos\left(N\theta\right)d\theta\\
 & = & \frac{2L+1}{2}\frac{i^{-N}}{\pi}\int_{0}^{\pi}\left[\left(\sum_{l'=0}^{\infty}(2l'+1)i^{l'}j_{l'}\left(k\cos\theta\right)\frac{2}{2l'+1}\delta_{l'L}\right)\right]\cos\left(N\theta\right)d\theta\\
 & = & \left(2L+1\right)\frac{i^{L-N}}{\pi}\int_{0}^{\pi}j_{L}\left(k\cos\theta\right)\cos\left(N\theta\right)d\theta\;.
\end{array}\label{eq:scriptj0}
\end{equation}

Using the series expansion \cite{wolfram.com/03.21.06.0019.01}

\begin{equation}
j_{l}(x)=\frac{1}{2}\sqrt{\pi}\left(\frac{x}{2}\right)^{l}\sum_{M=0}^{\infty}\frac{\left(-\frac{1}{4}\right)^{M}x^{2M}}{M!\Gamma\left(l+M+\frac{3}{2}\right)}\label{eq:wolfram03.21.06.0019.01}
\end{equation}
 this becomes

\begin{eqnarray}
a_{LN}\left(k\right) \hspace{-.3cm} & = & \hspace{-.3cm} \frac{(2L+1)}{\sqrt{\pi}}i^{L-\text{NN}}2^{-L-1}\left(\sum_{M=0}^{\text{\ensuremath{\infty}}}\frac{\left(-\frac{1}{4}\right)^{M}k^{L+2M}}{M!\Gamma\left(L+M+\frac{3}{2}\right)}\right)\int_{0}^{\pi}\cos^{L+2M}(\theta)\cos\left(\text{N}\theta\right)\, d\theta .\label{eq:a_w_M_sum}
\end{eqnarray}
 Gr\"{o}bner and Hofreiter\cite{GH p. 110 No. 332.14a} extended an
integral over the interval $\left[0,\frac{\pi}{2}\right]$ that has
three branches, to the interval $\left[0,\pi\right]$ with a prefactor
$(1+(-1)^{m+n})$ that renders the central one of the three possibilities
nonzero only for even values for $m+n$.

\begin{equation}
\int_{0}^{\pi}cos^{m}\theta\, cos(n\theta)d\theta=\left(1+(-1)^{m+n}\right)\frac{\pi}{2^{m+1}}\binom{m}{\frac{m-n}{2}}\quad\left[\mathrm{where\:}m\geq n>-1,\, m-n=2K\right]\;.\label{eq:GR5 p. 417 No. 3.631.17}
\end{equation}
 The other two branches, being for odd $m+n$ on $\left[0,\frac{\pi}{2}\right]$,
are zero on $\left[0,\pi\right]$ when this prefactor is included.
(Numerical integration also confirms that the contributions from the
$\left[\frac{\pi}{2},\pi\right]$ interval cancels the contributions
from the $\left[0,\frac{\pi}{2}\right]$ on these branches.) The fifth
edition of Gradshteyn and Ryzhik\cite{GR5 p. 417 No. 3.631.17} (in
which \emph{m} and \emph{n} are reversed in meaning) nevertheless
included all three branches and this prefactor on the interval $\left[0,\pi\right]$.
By their seventh edition, Gradshteyn and Ryzhik removed this integral
entirely despite the correctness of the central branch on the interval
$\left[0,\pi\right]$. Neither source noted the lower limit on \emph{m}
that we found, that $m\geq n>-1$.

The final form for the coefficient set of the Fourier-Legendre series
for the Bessel function $J_{N}\left(kx\right)$ is then 

\begin{equation}
\begin{array}{ccc}
a_{LN}\left(k\right) & = & \sqrt{\pi}(2L+1)2^{-L-1}i^{L-\text{N}}\sum_{M=0}^{\infty}\frac{\left(\left(-\frac{1}{4}\right)^{M}k^{L+2M}\right)}{2^{L+2M+1}\left(M!\Gamma\left(L+M+\frac{3}{2}\right)\right)}\\
 & \times & \left(1+(-1)^{L+2M+\text{N}}\right)\binom{L+2M}{\frac{1}{2}(L+2M-\text{N})}\\
 & = & \frac{\sqrt{\pi}2^{-2L-2}(2L+1)k^{L}i^{L-N}}{\Gamma\left(\frac{1}{2}(2L+3)\right)}\left(1+(-1)^{L+N}\right)\binom{L}{\frac{L-N}{2}}\\
 & \times & \,_{2}F_{3}\left(\frac{L}{2}+\frac{1}{2},\frac{L}{2}+1;L+\frac{3}{2},\frac{L}{2}-\frac{N}{2}+1,\frac{L}{2}+\frac{N}{2}+1;-\frac{k^{2}}{4}\right)\\
 & = & \sqrt{\pi}\,2^{-2L-2}(2L+1)k^{L}i^{L-N}\left(1+(-1)^{L+N}\right)\Gamma(L+1)\\
 & \times & \,_{2}\tilde{F}_{3}\left(\frac{L}{2}+\frac{1}{2},\frac{L}{2}+1;L+\frac{3}{2},\frac{L}{2}-\frac{N}{2}+1,\frac{L}{2}+\frac{N}{2}+1;-\frac{k^{2}}{4}\right)
\end{array}\;,\label{eq:a_L_as_2F3}
\end{equation}
where the final two steps are new with the present work. We have included
the final form using regularized hypergeometric functions \cite{wolfram.com/07.26.26.0001.01}

\begin{equation}
\,_{2}F_{3}\left(a_{1},a_{2};b_{1},b_{2},b_{3};z\right)=\Gamma\left(b_{1}\right)\Gamma\left(b_{2}\right)\Gamma\left(b_{3}\right)\,_{2}\tilde{F}_{3}\left(a_{1},a_{2};b_{1},b_{2},b_{3};z\right)\label{eq:regularized}
\end{equation}
and cancelled the $\Gamma\left(b_{i}\right)$ with gamma functions
in the denominators of the prefactors. Each give infinities that arise
whenever \emph{N\textgreater{}1} is an integer larger than \emph{L},
and of the same parity, resulting in indeterminacies in a computation
when one tries to use the conventional form of the hypergeometric
function.

For the special cases of $N=0,\,1$ the order of the hypergeometric
functions is reduced since the parameters $a_{2}=b_{3}$ and $a_{1}=b_{2}$,
resp., giving

\begin{equation}
\begin{array}{ccc}
a_{L0}\left(k\right) & = & \frac{\sqrt{\pi}i^{L}2^{-2L-2}(2L+1)k^{L}}{\Gamma\left(\frac{1}{2}(2L+3)\right)}\left(1+(-1)^{L}\right)\binom{L}{\frac{L}{2}}\\
 & \times & \,_{1}F_{2}\left(\frac{L}{2}+\frac{1}{2};\frac{L}{2}+1,L+\frac{3}{2};-\frac{k^{2}}{4}\right)\\
 & = & \sqrt{\pi}i^{L}2^{-2L-2}(2L+1)k^{L}\Gamma\left(\frac{L}{2}+1\right)\left(1+(-1)^{L}\right)\binom{L}{\frac{L}{2}}\\
 & \times & \,_{1}\tilde{F}_{2}\left(\frac{L}{2}+\frac{1}{2};\frac{L}{2}+1,L+\frac{3}{2};-\frac{k^{2}}{4}\right)
\end{array}\;,\label{eq:a_L0}
\end{equation}

and

\begin{equation}
\begin{array}{ccc}
a_{L1}\left(k\right) & = & \frac{\sqrt{\pi}i^{L-1}2^{-2L-2}(2L+1)k^{L}}{\Gamma\left(\frac{1}{2}(2L+3)\right)}\left(1+(-1)^{L+1}\right)\binom{L}{\frac{L-1}{2}}\\
 & \times & \,_{1}F_{2}\left(\frac{L}{2}+1;\frac{L}{2}+\frac{3}{2},L+\frac{3}{2};-\frac{k^{2}}{4}\right)\\
 & = & i^{L-1}2^{-L-2}(2L+1)k^{L}\Gamma\left(\frac{L}{2}+1\right)\left(1+(-1)^{L+1}\right)\\
 & \times & \,_{1}\tilde{F}_{2}\left(\frac{L}{2}+1;\frac{L}{2}+\frac{3}{2},L+\frac{3}{2};-\frac{k^{2}}{4}\right)
\end{array}\;,\label{eq:a_L1}
\end{equation}
In each special case the first form involving a hypergeometric function
has no indeterminacies, but we include the regularized hypergeometric
function version for completeness.

\section{A Numerical Check}

Numerical checks are an essential element of any analytical work,
particularly in checking convergence. Much of my research in physics
has revolved around writing my own programs to diagonalize very large
Hamiltonian matrices and I have found that I routinely need to calculate
in quadruple precision (33-digits) to gain consistent double-precision
results. Since this research often contains Bessel functions, it made
sense to craft a quadruple precision approximation for $J_{0}\left(kx\right)$that
would both allow for numerical checks and  be useful for some
future project of mine, or of other researchers, calculated in quadruple
precision. Though Castellanos and Rosenthal obtain fifteen-decimal
accuracy with their Rational Chebyshev Approximations,\cite{Rational Chebyshev}
and there is a ZBESJ Bessel function double precision package that
has 18-digit accuracy (UNIVAC double precision),\cite{Amos_double}
I found no prepackaged quadruple precision routines in a search of
software libraries. For this reason I have displayed all coefficients
in what follows in 33-digits, accuracy that matches IEEE selected\_real\_kind(33,
4931) extended precision%
\footnote{The second argument, 4931, indicates a range $10^{-4931}$to $10^{4931}-1$
using128 bits.\cite{floating_point}%
} available on some compilers. These results are given in the appendix
so that that readers may simply copy and paste them into a calculation routine.

The first 22 terms in the sum (\ref{eq:Fourier-Legendre series})
are then,  with $k=1$ (using $\text{E-}7$ as the programming shorthand for $\times10^{-7})$,

\begin{eqnarray}
J_{0}\left(x\right) & \cong & 0.9197304100897602393144211940806200P_{0}(x)\nonumber \\
 & - & 0.1579420586258518875737139671443637P_{2}(x)\nonumber \\
 & + & 0.003438400944601109232996887872072915P_{4}(x)\nonumber \\
 & - & 0.00002919721848828729693660590986125663P_{6}(x)\nonumber \\
 & + & 1.317356952447780977655616563143280\text{E-}7\: P_{8}(x)\nonumber \\
 & - & 3.684500844208203027173771096058866\text{E-}10\: P_{10}(x)\nonumber \\
 & + & 7.011830032993845928208803328211457\text{E-}13\; P_{12}(x)\nonumber \\
 & - & 9.665964369858912263671995372753346\text{E-}16\; P_{14}(x)\nonumber \\
 & + & 1.009636276824546446525342170924936\text{E-}18\; P_{16}(x)\nonumber \\
 & - & 8.266656955927637858991972584174117\text{E-}22\; P_{18}(x)\nonumber \\
 & + & 5.448244867762758725890082837839430\text{E-}25\; P_{20}(x)\nonumber \\
 & - & 2.952527182137354751675774606663400\text{E-}28\; P_{22}(x)\nonumber \\
 & + & 1.338856158858534469080898670096200\text{E-}31\; P_{24}(x)\nonumber \\
 & - & 5.154913186088512926193234837816582\text{E-}35\; P_{26}(x)\nonumber \\
 & + & 1.706231577038503450138564028467634\text{E-}38\; P_{28}(x)\nonumber \\
 & - & 4.906893556427796857473097979568289\text{E-}42\; P_{30}(x)\nonumber \\
 & + & 1.237489200717479383020539576221293\text{E-}45\; P_{32}(x)\nonumber \\
 & - & 2.759056237537871868604555688548364\text{E-}49\; P_{34}(x)\nonumber \\
 & + & 5.477382207172712629199714648396409\text{E-}53\; P_{36}(x)\nonumber \\
 & - & 9.744200345578852550688946057050674\text{E-}57\; P_{38}(x)\nonumber \\
 & + & 1.562280711659504489828025148995770\text{E-}60\; P_{40}(x)\nonumber \\
 & - & 2.269056283827394368836057470594599\text{E-}60\; P_{42}(x)\;.\label{eq:Fourier-Legendre series-numerical}
\end{eqnarray}

If we wish to check the convergence of this series at, say, the $\left|\varepsilon\right|<5\times10^{-8}$
level of the polynomial approximation given by E. E. Allen (and
reproduced in Abramowitz and Stegun),\cite{AS p. 369 No. 9.4.1} 
\begin{eqnarray}
J_{0}\left(x\right)& \cong &1-2.25\left(\frac{x}{3}\right)^{2}+1.26562\left(\frac{x}{3}\right)^{4}-0.316387\left(\frac{x}{3}\right)^{6}+0.0444479\left(\frac{x}{3}\right)^{8}\nonumber \\
&- &0.0039444\left(\frac{x}{3}\right)^{10} +0.0002100\left(\frac{x}{3}\right)^{12}\;,\label{eq:AS p. 369 No. 9.4.1}
\end{eqnarray}
at the latter's limiting range of $x=3$ (with a value of $-0.260052$),
we find that truncating the series after the $P_{12}(x)$ term (to
get contributions through $x^{12}$) is insufficient, giving two fewer
digits of accuracy ($-0.260045$). This shows that the optimization
Allen must have done to obtain his approximation has a significant
effect.

Testing accuracy at the 15-digit level shows that we can truncating
the series after the $P_{24}(x)$ term in the result, $-0.260051954901933$.
This $\left|\varepsilon\right|<5\times10^{-16}$ is much better than
one would expect from doubling the number of powers in our unoptimized approximation with $P_{12}(x)$.
In addition, the \emph{range} of applicability of the new Legendre
approximation truncated after $P_{24}(x)\sim x^{24}$, if one is satisfied
with $\left|\varepsilon\right|<5\times10^{-8}$, roughly doubles to
$J_{0}\left(6\right)\cong0.15064531$.

Finally, in our quality checking step, we found that if we imported
a 33-digit version of (\ref{eq:Fourier-Legendre series-numerical})
back into the software we had been using to generate it, \emph{Mathematica}
\emph{7}, we only obtained a result accurate to$\left|\varepsilon\right|<1\times10^{-32}$,
but if we imported the 34-digit version actually displayed
in (\ref{eq:Fourier-Legendre series-numerical}), we obtained a result
accurate to$\left|\varepsilon\right|<1\times10^{-34}$ over the range
$-3\leq x\leq3$, with $J_{0}\left(3\right)\cong-0.260051954901933437624154695977331$.
Accuracy exceeds $\left|\varepsilon\right|<7\times10^{-15}$ over
the range $-8\leq x\leq8,$ with $J_{0}\left(8\right)\cong0.171650807137560$.

\section{Summing a set of infinite series}

If  Allen's polynomial approximation \cite{AS p. 369 No. 9.4.1}
has its powers of three folded into the coefficients, it becomes
\begin{eqnarray}
J_{0}\left(x\right)\cong1-0.25\: x^{2}+0.0156249\, x^{4}-0.000434001\: x^{6}+6.77456\text{E-}6\, x^{8}-6.6799\text{E-}8\, x^{10}   \nonumber \\
+3.952\text{E-}10\, x^{12}    \;. \label{eq:AS p. 369 No. 9.4.1^x}
\end{eqnarray}
Noting that the third term is $1/64.0004$, one is led to wonder
if each of these coefficients is made up of inverse powers of primes
if one were to rederive a similar approximation using a process with
greater accuracy.

We can easily expand the Legendre polynomials into their constituent
terms and gather like powers in  (\ref{eq:Fourier-Legendre series-numerical})
to give an updated polynomial approximation,

\begin{eqnarray}
J_{0}\left(x\right) & \cong & 1.000000000000000000000000000000000000\, x^{0}\nonumber \\
 & - & 0.2500000000000000000000000000000000000\, x^{2}\nonumber \\
 & + & 0.01562500000000000000000000000000000000\, x^{4}\nonumber \\
 & - & 0.0004340277777777777777777777777777777778\, x^{6}\nonumber \\
 & + & 6.781684027777777777777777777777777778\text{E-}6\; x^{8}\nonumber \\
 & - & 6.781684027777777777777777777777777778\text{E-}8\; x^{10}\nonumber \\
 & + & 4.709502797067901234567901234567901235\text{E-}10\; x^{12}\nonumber \\
 & - & 2.402807549524439405391786344167296548\text{E-}12\; x^{14}\nonumber \\
 & + & 9.385966990329841427311665406903502142\text{E-}15\; x^{16}\nonumber \\
 & - & 2.896903392077111551639402903365278439\text{E-}17\; x^{18}\nonumber \\
 & + & 7.242258480192778879098507258413196097\text{E-}20\; x^{20}\nonumber \\
 & - & 1.496334396734045222954237036862230599\text{E-}22\; x^{22}\nonumber \\
 & + & 2.597802772107717400962217077885817011\text{E-}25\; x^{24}\nonumber \\
 & - & 3.842903509035084912666001594505646466\text{E-}28\, x^{26}\nonumber \\
 & + & 4.901662639075363409012757135849038860\text{E-}31\; x^{28}\nonumber \\
 & - & 5.446291821194848232236396817610043178\text{E-}34\; x^{30}\nonumber \\
 & + & 5.318644356635593976793356267197307791\text{E-}37\; x^{32}\nonumber \\
 & - & 4.600903422695150498956190542558224733\text{E-}40\; x^{34}\nonumber \\
 & + & 3.550079801462307483762492702591222788\text{E-}43\; x^{36}\nonumber \\
 & - & 2.458504017633176927813360597362342651\text{E-}46\; x^{38}\nonumber \\
 & + & 1.5365650110207355798833503733514641567\text{E-}49\; x^{40}\nonumber \\
 & - & 8.7106860035189091830121903251216788929\text{E-}53\; x^{42}\nonumber \\
 & \cong & 1-\frac{x^{2}}{2^{2}}+\frac{x^{4}}{2^{6}}-\frac{x^{6}}{2^{8}3^{2}}+\frac{x^{8}}{2^{14}3^{2}}-\frac{x^{10}}{2^{16}3^{2}5^{2}}+\frac{x^{12}}{2^{20}3^{4}5^{2}}-\frac{x^{14}}{2^{22}3^{4}5^{2}7^{2}}+\frac{x^{16}}{2^{30}3^{4}5^{2}7^{2}}\nonumber \\
 & - & \frac{x^{18}}{2^{32}3^{8}5^{2}7^{2}}+\frac{x^{20}}{2^{36}3^{8}5^{4}7^{2}}-\frac{x^{22}}{2^{38}3^{8}5^{4}7^{2}11^{2}}+\frac{x^{24}}{2^{44}3^{10}5^{4}7^{2}11^{2}}-\frac{x^{26}}{2^{46}3^{10}5^{4}7^{2}11^{2}13^{2}}\nonumber \\
 & + & \frac{x^{28}}{2^{50}3^{10}5^{4}7^{4}11^{2}13^{2}}-\frac{x^{30}}{2^{52}3^{12}5^{6}7^{4}11^{2}13^{2}}+\frac{x^{32}}{2^{62}3^{12}5^{6}7^{4}11^{2}13^{2}}\nonumber \\
 & - &\frac{x^{34}}{2^{64}3^{12}5^{6}7^{4}11^{2}13^{2}17^{2}}+  \frac{x^{36}}{2^{68}3^{16}5^{6}7^{4}11^{2}13^{2}17^{2}}-\frac{x^{38}}{2^{70}3^{16}5^{6}7^{4}11^{2}13^{2}17^{2}19^{2}}\nonumber \\
 & + & \frac{x^{40}}{2^{76}3^{16}5^{8}7^{4}11^{2}13^{2}17^{2}19^{2}}-\frac{x^{42}}{2^{78}3^{18}5^{8}7^{6}11^{2}13^{2}17^{2}19^{2}}\;.\label{eq:poly_approx_48_digit}
\end{eqnarray}
 One can see from the repeating digits that the fourth through sixth
lines in  (\ref{eq:poly_approx_48_digit}) have inverses that are
powers of primes, and one can even see that the fifth term is $2^{2}5^{2}$
times the sixth. Subsequent terms are not at all obviously inverse
powers of primes. Indeed, revealing those powers as analytic entities
required many terms in excess of those required to achieve quadruple
precision in the numerical results, and significantly higher precision.

One might suppose that that including series terms through $P_{24}(x)$
would be sufficient for revealing the inverse powers in the coefficient
of, say, the $x^{16}$ term, but that is not the case. Increasing the precision
from 33 to 48-digits did not improve the situation enough. At that
series truncation and with 48-digit precision, the inverse of the
coefficient of the $x^{16}$ term is $\mathbf{106542032486495}.616348409991752462411671619456197,$
whose integer part is in bold face. (We used the algebra and calculus
computer software \emph{Mathematica} 7 for this work.) This is not
a product of low-level primes. Adding one more term, $P_{26}(x)$,
is sufficient to bring it to $\mathbf{106542032486400}.113376300998684305400416345779209$,
whose integer part is $2^{30}3^{4}5^{2}7^{2}$. One additional term
makes this $\mathbf{106542032486400}.000104784167278249059923631013442$
and with every additional term added, the integer part remains the
same while the fractional part diminishes by several decimal places.
The coefficient of the $x^{24}$ term required 40 terms in the series,
with 48-digit precision, to establish convergence. The coefficient
of the $x^{42}$ term required 74 terms in the series, with 50-digit
precision, to establish convergence.

All of the coefficients in (\ref{eq:poly_approx_48_digit}) include
contributions from all 74 terms in the series, calculated with 50-digit
precision. These were then truncated to the 37-digit precision displayed
therein, except for the last two that required 38 digits. Taking their
reciprocals after truncation did not forestall revealing their integer
parts as powers of primes. The truncated power series gives $J_{0}\left(3\right)\cong-0.2600519549019334376241546959773314809$,
which matches \emph{Mathematica}'s BesselJ{[}0,3{]}=$-0.2600519549019334376241546959773314368$
(when set to 37-digit precision) within an error of $\left|\varepsilon\right|<4\times10^{-35}$
, as does the inverse prime version, 

\noindent
$J_{0}\left(3\right)\cong-\frac{90658024929169559805594876257679495662633}{348615048725002045174179287542005760000000}=-0.2600519549019334376241546959773314809$.

Close examination of the inverse prime version shows that 
we have managed to translate into integer powers the denominators
of the first 22 terms of the well-known series representation\cite{GR5 p. 970 No. 8.440}
\begin{equation}
J_{\nu}\left(x\right)=\sum_{k=0}^{\infty}\frac{(-1)^{k}\left(\frac{x}{2}\right)^{2k+\nu}}{k!\Gamma(k+\nu+1)}.\label{eq:GR5 p. 970 No. 8.440}
\end{equation}
(There is a fascinating analogue to this result arising from studies
of the Bessel difference equation.\cite{Cuchta}) But the outcome
worth the trouble of this investigation is that this process yields
a set of infinite sums whose values are inverse powers of primes.
The question of why so many Legendre polynomials are required to get
sufficiently accurate coefficients in (\ref{eq:poly_approx_48_digit})
led to an examination of precisely how a given $x^{n}$ contributor
from each Legendre polynomial contributes to the respective coefficient.

Looking back at the coefficients in  (\ref{eq:Fourier-Legendre series-numerical})
when multiplied by the constant terms in the Legendre polynomials
that multiply them, whose first few are 
\begin{eqnarray}
\left\{ P_{0}(x)=1,\; P_{2}(x)=\frac{1}{2}\left(3x^{2}-1\right),\; P_{4}(x)=\frac{1}{8}\left(35x^{4}-30x^{2}+3\right),\right.\nonumber \\
\left.\; P_{6}(x)=\frac{1}{16}\left(231x^{6}-315x^{4}+105x^{2}-5\right)\right\} \;,\label{eq:P0toP6}
\end{eqnarray}
 there is no reason to suspect that 
\begin{eqnarray}
0.919730410089760239314421194080620\nonumber \\
-\frac{1}{2}\left(-0.157942058625851887573713967144364\right)\nonumber \\
+\frac{3}{8}\left(0.00343840094460110923299688787207292\right)\nonumber \\
-\frac{5}{16}\left(-0.0000291972184882872969366059098612566\right)\nonumber \\
+\cdots & = & 1\label{eq:first4products}
\end{eqnarray}
 rather than some other number close to \emph{1}, but one sees uniform
convergence up through the accuracy of the calculation as one adds
additional terms, as seen in Table 1.

\vspace{.3cm}

\begin{tabular}{|c|}
\hline 
0.919730410089760239314421194080619970661964806513\tabularnewline
\hline 
\hline 
0.998701439402686183101278177652801821334364120020\tabularnewline
\hline 
0.999990839756911599063652010604829164640891568430\tabularnewline
\hline 
0.999999963887689188843944699951660807338623340119\tabularnewline
\hline 
0.999999999909168357337955807722426205787682516277\tabularnewline
\hline 
0.999999999999841620300892054094280728854756176645\tabularnewline
\hline 
0.99999999999999979732605041136082528291421094660\tabularnewline
\hline 
0.999999999999999999801573588581204243621535893434\tabularnewline
\hline 
0.999999999999999999999846581786952424267000275450\tabularnewline
\hline 
0.999999999999999999999999903953851991585994488442\tabularnewline
\hline 
0.999999999999999999999999999950320420042897103088\tabularnewline
\hline 
0.999999999999999999999999999999978412290228685090\tabularnewline
\hline 
0.999999999999999999999999999999999992008312509362\tabularnewline
\hline 
0.999999999999999999999999999999999999997449396440\tabularnewline
\hline 
0.999999999999999999999999999999999999999999290955\tabularnewline
\hline 
0.999999999999999999999999999999999999999999999827\tabularnewline
\hline 
1.00000000000000000000000000000000000000000000000\tabularnewline
\hline 
\end{tabular}

Table 1. The constant term of the Legendre series approximation as
increasing numbers of terms are added from  (\ref{eq:Fourier-Legendre series-numerical}),
to 48 digit accuracy.

\vspace{0.4cm}

One may more formally concluded that
\begin{equation}
\sum_{L=0}^{\infty}\,^{(2)}\frac{\sqrt{\pi}i^{L}(2L+1)\binom{L}{\frac{L}{2}}\binom{2L}{L}\left(\frac{1}{2}-\frac{L}{2}\right)_{\frac{L}{2}}\left(-\frac{L}{2}\right)_{\frac{L}{2}}\text{2}^{-3L-2}\,_{1}F_{2}\left(\frac{L}{2}+\frac{1}{2};\frac{L}{2}+1,L+\frac{3}{2};-\frac{1}{4}\right)}{\frac{L}{2}!\left(\frac{1}{2}-L\right)_{\frac{L}{2}}\Gamma\left(\frac{1}{2}(2L+3)\right)}=1\label{eq:1=00003D00003D}
\end{equation}
where the superscript ``(2)'' on the sum indicates one is summing
even values only (or one may retain the factor $\left(1+(-1)^{L}\right)$
in the sum as we do for similar equations below), which is a result we have not seen in the literature.
The Pochhammer symbols $\left(\frac{1}{2}-\frac{L}{2}\right)_{\frac{L}{2}}$
and so on derive from a shift to\cite{PBM3 p. 468 No. 7.3.1.206}
\begin{equation}
P_{n}(x)=2^{-n}\binom{2n}{n}x^{n}\,_{2}F_{1}\left(\frac{1}{2}-\frac{n}{2},-\frac{n}{2};\frac{1}{2}-n;\frac{1}{x^{2}}\right)\label{eq:PBM3 p. 468 No. 7.3.1.206}
\end{equation}
 in the explicit sum.\cite{PBM3 p. 430 No. 7.2.1.1} 

One may check that  the series (\ref{eq:1=00003D00003D}) converges
under the Cauchy criterion. That is, if we call a given term in (\ref{eq:1=00003D00003D})
$u_{L}$ and define the partial sums

\begin{equation}
S_{n}=\sum_{L=0}^{n}u_{L}\;,\label{eq:S_n}
\end{equation}
then the numerical series (\ref{eq:1=00003D00003D}) is convergent
if for each $\varepsilon>0$ there is a number $N(\varepsilon)$
such that

\begin{equation}
\left|S_{m}-S_{n}\right|<\varepsilon\label{eq:Cauchy}
\end{equation}
for all $m>n>N$. If we take, say, $\varepsilon=10^{-21}$, we see
that the difference of the ninth and eighth lines $S_{10}-S_{9}$
in Table 1, above, fulfills this bound  (with $m=10$ and $n=9$). If we instead take $\varepsilon=2\times10^{-46}$,
we see that the difference of the last two lines $S_{16}-S_{15}$
fulfills this tighter bound. 

The summed series (\ref{eq:1=00003D00003D}) was a consequence of
our expansion of $J_{N}\left(kx\right)$ in a Fourier-Legendre series
after setting $k=1$. Including \emph{k} poses no problem despite
its appearance as the argumemt of the $\,_{1}F_{2}\left(\frac{L}{2}+\frac{1}{2};\frac{L}{2}+1,L+\frac{3}{2};-\frac{k^{2}}{4}\right)$
function as well as there being a $k^{L}$ factor in the argument
of the sum. It ends up contributing a very clean factor of $k^{2h}$
to the right-hand side, below. The 43 verified summed series given
by the present approach (with $0\leq h\leq42$) are

\begin{eqnarray}
\sum_{L=0}^{\infty}\, &  & \hspace{-0.9cm}\frac{\sqrt{\pi}i^{L}2^{-3L-2}\left(1+(-1)^{L}\right)(2L+1)\binom{L}{\frac{L}{2}}\binom{2L}{L}\,_{1}F_{2}\left(\frac{L}{2}+\frac{1}{2};\frac{L}{2}+1,L+\frac{3}{2};-\frac{k^{2}}{4}\right)}{\Gamma\left(\frac{1}{2}(2L+3)\right)\left(\frac{L}{2}-h\right)!}k^{L}\nonumber \\
 & \times & \left[\frac{\left(\frac{1}{2}-\frac{L}{2}\right)_{\frac{L}{2}-h}\left(-\frac{L}{2}\right)_{\frac{L}{2}-h}}{\left(\frac{1}{2}-L\right)_{\frac{L}{2}-h}}\right]= \frac{(-1)^{h}2^{-2h}}{h!\Gamma(h+1)}k^{2h}\;.\label{eq:1/primes^n=00003D00003D}
\end{eqnarray}

To verify the final, $h=42$, summed series, we had to take the upper
limit on the number of terms in the series $\geq h+74$ in order to
obtain a percent difference between left- and right-hand sides that
was $\leq10^{-33}$ because the first \emph{h} terms in the series
do not contribute. For $h=0$, an upper limit on the number of terms
in the series $\geq h+44$ was sufficient.

We wish to provide two alternative forms for readers seeking to sum
a series involving 

\noindent
$\,_{1}F_{2}\left(\frac{L}{2}+\frac{1}{2};\frac{L}{2}+1,L+\frac{3}{2};-\frac{k^{2}}{4}\right)$
hypergeometric functions, who may have somewhat different coefficients
than the above: One may use the relation\cite{dlmf 5.2.6} 
\begin{equation}
(-a)_{n}=(-1)^{n}(a-n+1)_{n}\label{eq:dlmf 5.2.6}
\end{equation}
 and the primary definition of the Pochhammer symbol\cite{http://functions.wolfram.com/06.10.02.0001.01}
\begin{equation}
(a)_{n}=\frac{\Gamma(a+n)}{\Gamma(a)}\text{/;}\neg(-a\in\mathbb{Z}\land-a\geq0\land n\in\mathbb{Z}\land n\leq-a)
\end{equation}
 to rewrite the term in square brackets in infinite sum (\ref{eq:1/primes^n=00003D00003D})
to give two alternative forms:

\begin{eqnarray}
\sum_{L=0}^{\infty}\, &  & \hspace{-0.9cm}\frac{\sqrt{\pi}i^{L}2^{-3L-2}\left(1+(-1)^{L}\right)(2L+1)\binom{L}{\frac{L}{2}}\binom{2L}{L}\,_{1}F_{2}\left(\frac{L}{2}+\frac{1}{2};\frac{L}{2}+1,L+\frac{3}{2};-\frac{k^{2}}{4}\right)}{\Gamma\left(\frac{1}{2}(2L+3)\right)\left(\frac{L}{2}-h\right)!}k^{L}\nonumber \\
 & \times & \left[\frac{(-1)^{\frac{L}{2}-h}\left(h+\frac{1}{2}\right)_{\frac{L}{2}-h}(h+1)_{\frac{L}{2}-h}}{\left(h+\frac{L}{2}+\frac{1}{2}\right)_{\frac{L}{2}-h}}\right]\nonumber \\
 & = & \sum_{L=0}^{\infty}\,\frac{\sqrt{\pi}i^{L}2^{-3L-2}\left(1+(-1)^{L}\right)(2L+1)\binom{L}{\frac{L}{2}}\binom{2L}{L}\,_{1}F_{2}\left(\frac{L}{2}+\frac{1}{2};\frac{L}{2}+1,L+\frac{3}{2};-\frac{k^{2}}{4}\right)}{\Gamma\left(\frac{1}{2}(2L+3)\right)\left(\frac{L}{2}-h\right)!}k^{L}\nonumber \\
 & \times & \left[\frac{(-1)^{\frac{L}{2}-h}2^{2h-L}\Gamma(L+1)\Gamma\left(h+\frac{L}{2}+\frac{1}{2}\right)}{\Gamma(2h+1)\Gamma\left(L+\frac{1}{2}\right)}\right]=\frac{(-1)^{h}2^{-2h}}{h!\Gamma(h+1)}k^{2h}\;.\label{eq:1/primes^n_gamma}
\end{eqnarray}

\section{Series arising from the $J_{1}\left(x\right)$ Fourier-Legendre series}

The first 22 terms in the $J_{1}\left(x\right)$ Fourier-Legendre
series (\ref{eq:Fourier-Legendre series}) are

\begin{eqnarray}
J_{1}\left(x\right) & \cong & 0.4635981705953810635941110039338702P_{1}(x)\nonumber \\
 & - & 0.02386534565840739796307209416484866P_{3}(x)\nonumber \\
 & + & 0.0003197243559720047638524757623256028P_{5}(x)\nonumber \\
 & - & 1.970519180666594250258062929391112\text{E-}6\; P_{7}(x)\nonumber \\
 & + & 6.987247473097807218791759410157014\text{E-}9\; P_{9}(x)\nonumber \\
 & - & 1.610500056046875027807002442953327\text{E-}11\; P_{11}(x)\nonumber \\
 & + & 2.607086592441628842939248193619909\text{E-}14\; P_{13}(x)\nonumber \\
 & - & 3.127311482540796882144713619567442\text{E-}17\; P_{15}(x)\nonumber \\
 & + & 2.891424081787050739827382596616064\text{E-}20\; P_{17}(x)\nonumber \\
 & - & 2.123664534779369199214414455720317\text{E-}23\; P_{19}(x)\nonumber \\
 & + & 1.269011201758673511714553707528186\text{E-}26\; P_{21}(x)\nonumber \\
 & - & 6.290201939135925763576871358738600\text{E-}30\; P_{23}(x)\nonumber \\
 & + & 2.628135796989325452573870774267213\text{E-}33\; P_{25}(x)\nonumber \\
 & - & 9.381575562723076109283258050667642\text{E-}37\; P_{27}(x)\nonumber \\
 & + & 2.894337242415984040941859061022419\text{E-}40\; P_{29}(x)\nonumber \\
 & - & 7.794444104104171684395094261174814\text{E-}44\; P_{31}(x)\nonumber \\
 & + & 1.848200759818170134895867052306767\text{E-}47\; P_{33}(x)\nonumber \\
 & - & 3.888249639773912225694535890329244\text{E-}51\; P_{35}(x)\nonumber \\
 & + & 7.306978718807123633044120058516188\text{E-}55\; P_{37}(x)\nonumber \\
 & - & 1.234022530456621571127590099647796\text{E-}58\; P_{39}(x)\nonumber \\
 & + & 1.883067799255568915649461884255428\text{E-}62\; P_{41}(x)\nonumber \\
 & - & 2.609122884536350861268195351045890\text{E-}66\; P_{43}(x)\label{eq:Fourier-Legendre series-numerical_J1}
\end{eqnarray}

If we wish to check the convergence of this series at, say, the $\left|\varepsilon\right|<1.3\times10^{-8}$
level given by the polynomial approximation given by E. E. Allen (and
reproduced in Abramowitz and Stegun),\cite{AS p. 370 No. 9.4.4},

\begin{eqnarray}
J_{1}\left(x\right) & \cong & \frac{x}{2}-0.062499983x^{3}+0.0026041448x^{5}-0.00005424265x^{7}+6.7568816\text{E-}7x^{9}\nonumber \\
 & - & 5.3788\text{E-}9x^{11}+2.087\text{E-}11x^{13}\;,\label{eq:AS  p. 370 No. 9.4.4}
\end{eqnarray}
at the latter's limiting range of $x=3$ (with a value of $0.339059$)
, we find that truncating the series after the $P_{13}(x)$ term (to
get contributions through $x^{13}$) is insufficient, giving one less
digit of accuracy ($0.339060$), again showing that the optimization
Allen must have done to obtain his approximation had a significant
effect.

Testing accuracy at the 15-digit level shows that we can truncating
the series after the $P_{25}(x)$ term in the result, $0.339058958525937$.
This $\left|\varepsilon\right|<4\times10^{-17}$ is  better than
one would expect for doubling the number of powers from the $P_{13}(x)$ truncation.
In addition, the \emph{range} of applicability of the new Legendre
polynomial approximation truncated after $P_{24}(x)\sim x^{24}$,
if one is satisfied with $\left|\varepsilon\right|<1\times10^{-8}$,
roughly doubles to $J_{0}\left(6\right)\cong-0.276684$.

Finally, we confirm that the series through $P_{43}(x)$ does give
33-digit accuracy for $-3\leq x\leq3$, with $J_{1}\left(3\right)\cong0.33905895852593645892551459720648$.
Accuracy exceeds $\left|\varepsilon\right|<9\times10^{-15}$ over
the range $-8\leq x\leq8,$ with $J_{1}\left(8\right)\cong0.234636346853916$.

We can again expand the Legendre polynomials into their constituent
terms and gather like powers in  (\ref{eq:Fourier-Legendre series-numerical_J1})
to give an updated polynomial approximation,

\begin{eqnarray}
J_{1}\left(x\right) & \cong & 0.5000000000000000000000000000000000000\, x\nonumber \\
 & - & 0.06250000000000000000000000000000000000\, x^{3}\nonumber \\
 & + & 0.002604166666666666666666666666666666667\, x^{5}\nonumber \\
 & - & 0.00005425347222222222222222222222222222222x^{7}\nonumber \\
 & + & 6.781684027777777777777777777777777778\text{E-}7\; x^{9}\nonumber \\
 & - & 5.651403356481481481481481481481481481\text{E-}9\; x^{11}\nonumber \\
 & + & 3.363930569334215167548500881834215168\text{E-}11\; x^{13}\nonumber \\
 & - & 1.501754718452774628369866465104560343\text{E-}13\; x^{15}\nonumber \\
 & + & 5.214426105738800792950925226057501190\text{E-}16\; x^{17}\nonumber \\
 & - & 1.448451696038555775819701451682639219\text{E-}18\; x^{19}\nonumber \\
 & + & 3.291935672814899490499321481096907317\text{E-}21\; x^{21}\nonumber \\
 & - & 6.234726653058521762309320986925960827\text{E-}24\; x^{23}\nonumber \\
 & + & 9.991549123491220772931604145714680813\text{E-}27\; x^{25}\nonumber \\
 & - & 1.372465538941101754523571998037730881\text{E-}29\; x^{27}\nonumber \\
 & + & 1.633887546358454469670919045283012953\text{E-}32\; x^{29}\nonumber \\
 & - & 1.701966194123390072573874005503138493\text{E-}35\; x^{31}\nonumber \\
 & + & 1.564307163716351169645104784469796409\text{E-}38\; x^{33}\nonumber \\
 & - & 1.278028728526430694154497372932840204\text{E-}41\; x^{35}\nonumber \\
 & + & 9.342315267006072325690770269976902073\text{E-}45\; x^{37}\nonumber \\
 & - & 6.146260044082942319533401493405856627\text{E-}48\; x^{39}\nonumber \\
 & + & 3.658488121477941856865119936551105135\text{E-}51\; x^{41}\nonumber \\
 & - & 1.979701364436115723411861437527654294\text{E-}54\; x^{43}\nonumber \\
 & = & \frac{x}{2}-\frac{x^{3}}{2^{4}}+\frac{x^{5}}{2^{7}3}-\frac{x^{7}}{2^{11}3^{2}}+\frac{x^{9}}{2^{15}3^{2}5}-\frac{x^{11}}{2^{18}3^{3}5^{2}}+\frac{x^{13}}{2^{21}3^{4}5^{2}7}-\frac{x^{15}}{2^{26}3^{4}5^{2}7^{2}}\nonumber \\
 & + & \frac{x^{17}}{2^{31}3^{6}5^{2}7^{2}}-\frac{x^{19}}{2^{34}3^{8}5^{3}7^{2}}+\frac{x^{21}}{2^{37}3^{8}5^{4}7^{2}11}-\frac{x^{23}}{2^{41}3^{9}5^{4}7^{2}11^{2}}\nonumber \\
 & + & \frac{x^{25}}{2^{45}3^{10}5^{4}7^{2}11^{2}13}-\frac{x^{27}}{2^{48}3^{10}5^{4}7^{3}11^{2}13^{2}}+\frac{x^{29}}{2^{51}3^{11}5^{5}7^{4}11^{2}13^{2}}-\frac{x^{31}}{2^{57}3^{12}5^{6}7^{4}11^{2}13^{2}}\nonumber \\
 & + & \frac{x^{33}}{2^{63}3^{12}5^{6}7^{4}11^{2}13^{2}17}-\frac{x^{35}}{2^{66}3^{14}5^{6}7^{4}11^{2}13^{2}17^{2}}+\frac{x^{37}}{2^{69}3^{16}5^{6}7^{4}11^{2}13^{2}17^{2}19} \;,\label{eq:J1series}  \\
 & - &\frac{x^{39}}{2^{73}3^{16}5^{7}7^{4}11^{2}13^{2}17^{2}19^{2}}+\frac{x^{41}}{2^{77}3^{17}5^{8}7^{5}11^{2}13^{2}17^{2}19^{2}}-\frac{x^{43}}{2^{80}3^{18}5^{8}7^{6}11^{3}13^{2}17^{2}19^{2}}\nonumber
\end{eqnarray}
with the latter form an integer-power realization of first 22 terms
of the well-known series representation\cite{GR5 p. 970 No. 8.440}.
This time all terms in the first half were truncated to the 37-digit
precision displayed therein. 

The truncated power series gives $J_{1}\left(3\right)\cong0.339058958525936458925514597206478894$,
which matches \emph{Mathematica}'s BesselJ{[}1,3{]}=$0.3390589585259364589255145972064788970$
(when set to 37-digit precision) within an error of $\left|\varepsilon\right|<3\times10^{-36}$,
as does the inverse prime version 

\noindent
$J_{1}\left(3\right)\cong\frac{23266944578863553712347684324898325104584007}{68622120117447770997443712389847449600000000}={\scriptstyle 0.3390589585259364589255145972064788941}$.

The Legendre series, via the above development , gives us another
set of infinite series (of which we have confirmed $0\leq h\leq43$),

\begin{eqnarray}
\sum_{L=1}^{\infty}\, &  & \hspace{-0.9cm}\frac{\sqrt{\pi}i^{L-1}\left(1+(-1)^{L+1}\right)(2L+1)2^{-3L-2}\binom{L}{\frac{L-1}{2}}\binom{2L}{L}(-1)^{-h+\frac{L}{2}-\frac{1}{2}}}{\Gamma\left(\frac{1}{2}(2L+3)\right)\left(-h+\frac{L}{2}-\frac{1}{2}\right)!} k^L \nonumber \\
 & \times &\,_{1}F_{2}\left(\frac{L}{2}+1;\frac{L}{2}+\frac{3}{2},L+\frac{3}{2};-\frac{k^2}{4}\right) \left[\frac{\left(\frac{1}{2}-\frac{L}{2}\right)_{\frac{L-1}{2}-h}\left(-\frac{L}{2}\right)_{\frac{L-1}{2}-h}}{\left(\frac{1}{2}-L\right)_{\frac{L-1}{2}-h}}\right]\nonumber \\
 & = & \sum_{L=1}^{\infty}\,^{(2)}\frac{\sqrt{\pi}i^{L-1}\left(1+(-1)^{L+1}\right)(2L+1)2^{-3L-2}\binom{L}{\frac{L-1}{2}}\binom{2L}{L}}{\Gamma\left(\frac{1}{2}(2L+3)\right)\left(-h+\frac{L}{2}-\frac{1}{2}\right)!}  k^L \nonumber \\
 & \times & \,_{1}F_{2}\left(\frac{L}{2}+1;\frac{L}{2}+\frac{3}{2},L+\frac{3}{2};-\frac{k^2}{4}\right) \left[(-1)^{-h+\frac{L}{2}-\frac{1}{2}}\frac{(h+1)_{-h+\frac{L}{2}-\frac{1}{2}}\left(h+\frac{3}{2}\right)_{-h+\frac{L}{2}-\frac{1}{2}}}{\left(h+\frac{L}{2}+1\right)_{-h+\frac{L}{2}-\frac{1}{2}}}\right]\nonumber \\
 & = & \sum_{L=1}^{\infty}\,^{(2)}\frac{\sqrt{\pi}i^{L-1}\left(1+(-1)^{L+1}\right)(2L+1)2^{-3L-2}\binom{L}{\frac{L-1}{2}}\binom{2L}{L}}{\Gamma\left(\frac{1}{2}(2L+3)\right)\left(-h+\frac{L}{2}-\frac{1}{2}\right)!}  k^L \nonumber \\
 & \times & \,_{1}F_{2}\left(\frac{L}{2}+1;\frac{L}{2}+\frac{3}{2},L+\frac{3}{2};-\frac{k^2}{4}\right) \left[(-1)^{-h+\frac{L}{2}-\frac{1}{2}}\frac{2^{2h-L+1}\Gamma(L+1)\Gamma\left(h+\frac{L}{2}+1\right)}{\Gamma(2h+2)\Gamma\left(L+\frac{1}{2}\right)}\right] \nonumber \\
 & = & \frac{(-1)^{h}2^{-2h-1}}{h!\Gamma(h+2)} k^{2 h +1}\;.\label{eq:J1summedtoprimes}
\end{eqnarray}

To verify the final summed series, with $h=43$, we had to take the upper
limit on the number of terms in the series $\geq h+78$ in order to
obtain a percent difference between left- and right-hand sides that
was $\leq10^{-33}$ because the first \emph{h} terms in the series
do not contribute. For $h=0$, an upper limit on the number of terms
in the series $\geq h+45$ was sufficient.

\section{Series arising from the $I_{n}\left(x\right)$ Fourier-Legendre series}

Because the modified Bessel functions of the first kind $I_{N}\left(kx\right)$
are related to the ordinary Bessel functions by the relation \cite{GR5 p. 961 No. 8.406.3}
\begin{equation}
I_{n}(z)=i^{-n}J_{n}(iz)\;,\label{eq:I_from_J}
\end{equation}
 we merely need to multiply by $i^{-n}$ and set $k=i$ in  (\ref{eq:a_L_as_2F3})
to obtain the $I_{0}\left(x\right)$ Fourier-Legendre series, the
first 24 terms of which are

\begin{eqnarray}
I_{0}\left(x\right) & \cong & 1.086521097023589815837941923492506P_{0}(x)\nonumber \\
 & + & 0.1758046819215242662605951354261250P_{2}(x)\nonumber \\
 & + & 0.003709009244052882533923838165527033P_{4}(x)\nonumber \\
 & + & 0.00003095105270992432198613744608777602\: P_{6}(x)\nonumber \\
 & + & 1.381259734719773538320052305224506\text{E-}7\: P_{8}(x)\nonumber \\
 & + & 3.834312601086373005317788906125573\text{E-}10\: P_{10}(x)\nonumber \\
 & + & 7.257172450096213936720667660411978\text{E-}13\; P_{12}(x)\nonumber \\
 & + & 9.962746978836018020128433111635975\text{E-}16\; P_{14}(x)\nonumber \\
 & + & 1.037251346110052630963705477046736\text{E-}18\; P_{16}(x)\nonumber \\
 & + & 8.470496863240475339343499321604116\text{E-}22\; P_{18}(x)\nonumber \\
 & + & 5.570541399858852219278260483523687\text{E-}25\; P_{20}(x)\nonumber \\
 & + & 3.013347383234528850224689041823201\text{E-}28\; P_{22}(x)\nonumber \\
 & + & 1.364338005353527272638479093175249\text{E-}31\; P_{24}(x)\nonumber \\
 & + & 5.246088467162281648944660989359565\text{E-}35\; P_{26}(x)\nonumber \\
 & + & 1.734417052236546525979562610169336\text{E-}38\; P_{28}(x)\nonumber \\
 & + & 4.982929889631203560686967762401821\text{E-}42\; P_{30}(x)\nonumber \\
 & + & 1.255546430559877621587201790700357\text{E-}45\; P_{32}(x)\nonumber \\
 & + & 2.797096596401706413444068821508193\text{E-}49P_{34}(x)\nonumber \\
 & + & 5.548955677049963483909673845071489\text{E-}53\; P_{36}(x)\nonumber \\
 & + & 9.865206225205083247212985096531573\text{E-}57\; P_{38}(x)\nonumber \\
 & + & 1.580763691652306983099443761944673\text{E-}60\; P_{40}(x)\nonumber \\
 & + & 2.294688331479205281600814719914093\text{E-}64\; P_{42}(x)\nonumber \\
 & + & 3.031771495580703895127109933386607\text{E-}68\; P_{44}(x)\nonumber \\
 & + & 3.661200772680598752990852186025167\text{E-}72\; P_{46}(x)\;.\label{eq:Fourier-Legendre seriesI0}
\end{eqnarray}

If we wish to check the convergence of this series at, say, the $\left|\varepsilon\right|<1.6\times10^{-7}$
level of the polynomial approximation given by E. E. Allen (and
reproduced in Abramowitz and Stegun) \cite{AS p. 378 No. 9.8.1},

\begin{eqnarray}
I_{0}\left(x\right) & \cong & 1+3.5156229\left(\frac{x}{3.75}\right)^{2}+3.0899424\left(\frac{x}{3.75}\right)^{4}+1.2067492\left(\frac{x}{3.75}\right)^{6}\nonumber \\
 & + & 0.2659732\left(\frac{x}{3.75}\right)^{8}+0.0360768\left(\frac{x}{3.75}\right)^{10}+0.0045813\left(\frac{x}{3.75}\right)^{12}\label{eq:AS p. 378 No. 9.8.1}\\
 & = & 1+0.25\: x^{2}+0.0156252\: x^{4}+0.00043394\: x^{6}+6.801234\text{E-}6\; x^{8}+\text{6.56017\text{E-}8\;\ \ensuremath{x^{10}}}\nonumber \\
 & + & 5.9240\text{E-}10\; x^{12}\nonumber 
\end{eqnarray}

\noindent at the latter's limiting range of $x=3.75$ (with a value
of $9.11895$), we find that truncating the series after the $P_{12}(x)$
term (to get contributions through $x^{12}$) is insufficient, giving
two fewer digits of accuracy ($9.1187$). This shows that the optimization
Allen must have done to obtain his approximation has a significant
effect.

Testing accuracy at the 15-digit level shows that we can truncating
the series after the $P_{28}(x)$ term in the result, $9.1189458608445666$,
with $\left|\varepsilon\right|<5\times10^{-17}$. In addition, the
\emph{range} of applicability of the new Legendre approximation truncated
after $P_{24}(x)\sim x^{24}$, if one is satisfied with $\left|\varepsilon\right|<1.6\times10^{-7}$
, roughly doubles to $I_{0}\left(7.5\right)\cong268.1613$. 

Finally, in our quality checking step, we found that if we imported
a 33-digit version of (\ref{eq:Fourier-Legendre seriesI0}) back into
the software we had been using to generate it, \emph{Mathematica}
\emph{7}, we only obtained a result accurate to$\left|\varepsilon\right|<1\times10^{-32}$,
but if we imported the 34-digit version version actually displayed
in (\ref{eq:Fourier-Legendre seriesI0}), we obtained a result accurate
to $\left|\varepsilon\right|<1\times10^{-34}$ over the range $-3.75\leq x\leq3.75$,
with $I_{0}\left(3.75\right)\cong9.118945860844566690670997606599715$.
Accuracy exceeds $\left|\varepsilon\right|<6\times10^{-14}$ over
the range $-8\leq x\leq8,$ with $I_{0}\left(8\right)\cong427.5641157218$.

If we turn our attention to the latter form of the polynomial approximation
(\ref{eq:AS p. 378 No. 9.8.1}), and compare it with the $J_{0}$
polynomial approximation (\ref{eq:AS p. 369 No. 9.4.1^x}) one sees
that the  $I_{0}$ version is approximately the $J_{0}$ version
with all of the negative signs reversed. That the correspondence is
not exact for the higher-power terms likely is a result of the optimization
scheme in the two cases having slightly different ranges of validity.

Indeed, if we apply (\ref{eq:I_from_J}) to (\ref{eq:GR5 p. 970 No. 8.440}), each numerator becomes $i^{\nu}(-1)^{k}i^{2k+\nu}\left(\frac{x}{2}\right)^{2k+\nu}=(-1)^{\nu}\left(\frac{x}{2}\right)^{2k+\nu}$which
 indeed gives the series \cite{GR5 p. 971 No. 8.445} 

\begin{equation}
I_{\nu}\left(x\right)=\sum_{k=0}^{\infty}\frac{\left(\frac{x}{2}\right)^{2k+\nu}}{k!\Gamma(k+\nu+1)}\;.\label{eq:GR5 p. 971 No. 8.445}
\end{equation}
We, thus, need not display our 24-term polynomial approximation for
$I_{0}$ apart from the two terms unneeded for the  the $J_{0}$ version
that are necessary for   $I_{0}$:

\begin{eqnarray}
... &+& 4.4993212828093539168451396307446688497\text{E-}56 x^{44} \nonumber \\
& + &
 2.1263333094562164068266255343783879252\text{E-}59 x^{46}\;.\label{eq:x^44,46}
\end{eqnarray}

\noindent
One obtains the other terms by simply negating the negative signs in
(\ref{eq:poly_approx_48_digit}). This is not true of (\ref{eq:Fourier-Legendre seriesI0})
because the arguments of the Legendre polynomials do not undergo $x\rightarrow ix$
since they derive from the definition of the Fourier-Legendre series,
(\ref{eq:Fourier-Legendre series}). The \emph{k}-dependence is entirely
within the coefficients $a_{LN}\left(k\right)$. 

Furthermore, the $I_{0}$ Legendre series expansion leads to no new
set of summed series since these would simply be (\ref{eq:1/primes^n=00003D00003D})
with $k=i\kappa$. 

Although  (\ref{eq:Fourier-Legendre series}) allows one to easily
compute the Fourier-Legendre series for any $J_{n}\left(x\right)$
or $I_{n}\left(x\right)$, to enable readers to find these series
for higher indices by recursion \cite{GR5 p. 979 No. 8.471.1 & p. 981 No. 8.486.1}
we give the first 24 terms in the $I_{1}\left(x\right)$ Fourier-Legendre
series to complete the required pair:%

\begin{eqnarray}
I_{1}\left(x\right) & \cong & 0.5386343421852555592809081051666336\, P_{1}(x)\nonumber \\
 & + & 0.02618069164825977449795296407260333\, P_{3}(x)\nonumber \\
 & + & 0.0003419851912550806236210094361507344\, P_{5}(x)\nonumber \\
 & + & 2.077651971699656963860267070724864\, \text{E-}6\; P_{7}(x)\nonumber \\
 & + & 7.299001518662431414905576324932877\, \text{E-}9\; P_{9}(x)\nonumber \\
 & + & 1.671443482954853739162527767203215\, \text{E-}11\; P_{11}(x)\nonumber \\
 & + & 2.692744551459235232734936666452704\, \text{E-}14\; P_{13}(x)\nonumber \\
 & + & 3.218106754771162455853759838545282\, \text{E-}17\; P_{15}(x)\nonumber \\
 & + & 2.966624646773403074824196435937542\, \text{E-}20\; P_{17}(x)\nonumber \\
 & + & 2.173686883720901031568436047655748\, \text{E-}23\; P_{19}(x)\nonumber \\
 & + & 1.296326265789554875546711294998344\, \text{E-}26\; P_{21}(x)\nonumber \\
 & + & 6.414855102151415733596588296578833\, \text{E-}30\; P_{23}(x)\nonumber \\
 & + & 2.676389925142875786863074285387196\, \text{E-}33\; P_{25}(x)\nonumber \\
 & + & 9.542035089444710700263714658817980\, \text{E-}37\; P_{27}(x)\nonumber \\
 & + & 2.940669572337884276201779460203862\, \text{E-}40\; P_{29}(x)\nonumber \\
 & + & 7.911705033029330504663434638574930\, \text{E-}44\; P_{31}(x)\nonumber \\
 & + & 1.874426565726980007813411585706840\, \text{E-}47\; P_{33}(x)\nonumber \\
 & + & 3.940459072597980181771454632976250\, \text{E-}51\; P_{35}(x)\nonumber \\
 & + & 7.400090413796917559009186360838868\, \text{E-}55\; P_{37}(x)\nonumber \\
 & + & 1.248984620737396858084740490332061\, \text{E-}58\; P_{39}(x)\nonumber \\
 & + & 1.904842982553207494042180613785837\, \text{E-}62\; P_{41}(x)\nonumber \\
 & + & 2.637959760920312924684635466402215\, \text{E-}66\; P_{43}(x)\nonumber \\
 & + & 3.332061910821697596383220274010501\, \text{E-}70\; P_{45}(x)\:.\label{eq:Fourier-Legendre series-I1}
\end{eqnarray}

If we wish to check the convergence of this series at, say, the $\left|\varepsilon\right|<8\times10^{-9}$
level given by the polynomial approximation given by E. E. Allen (and
reproduced in Abramowitz and Stegun) \cite{AS p. 378 No. 9.8.3},

\begin{eqnarray}
I_{1}\left(x\right) & \cong & \frac{x}{2}+0.0625x^{3}+0.00260419x^{5}+0.0000542445x^{7}+6.79868\text{E-}7x^{9}+5.48303\text{E-}9x^{11}\nonumber \\
 & + & 4.191\text{E}-11x^{13}\;,\label{eq:AS p. 378 No. 9.8.3}
\end{eqnarray}

\noindent at the latter's limiting range of $x=3.75$ (with a value
of $7.78002$), we find that truncating the series after the $P_{13}(x)$
term (to get contributions through $x^{13}$) is insufficient, giving
three fewer digits of accuracy ($7.77996$). This again shows that the optimization
Allen must have done to obtain his approximation has a significant
effect.

Testing accuracy at the 15-digit level shows that we can truncating
the series after the $P_{27}(x)$ term in the result, $7.780015229824415$,
with $\left|\varepsilon\right|<4\times10^{-16}$. In addition, the
\emph{range} of applicability of the new Legendre approximation truncated
after $P_{27}(x)\sim x^{27}$, if one is satisfied with $\left|\varepsilon\right|<5\times10^{-9}$
, roughly doubles to $I_{1}\left(6.5\right)\cong97.735011$. 

Finally, in our quality checking step, we found that if we imported
a 33-digit version of (\ref{eq:Fourier-Legendre series-I1}) back into
the software we had been using to generate it, \emph{Mathematica}
\emph{7}, we only obtained a result accurate to$\left|\varepsilon\right|<1\times10^{-32}$,
but if we imported the 34-digit version  actually displayed
in (\ref{eq:Fourier-Legendre seriesI0}), we obtained a result accurate
to $\left|\varepsilon\right|<1\times10^{-34}$ over the range $-3.75\leq x\leq3.75$,
with $I_{1}\left(3.75\right)\cong7.780015229824415864988676277516113$.
Accuracy exceeds $\left|\varepsilon\right|<1\times10^{-18}$ over
the range $-8\leq x\leq8,$ with $I_{1}\left(8\right)\cong399.873136782560098$.

Because the above-noted correspondence between the power-series versions
of $J_{0}(x)$ and $I_{0}(x)$ (reversing all of the negative signs
in the former to achieve the latter) applies as well to $J_{1}(x)$
and $I_{1}(x)$, there is no need to display the power-series version
of the latter except for the additional term,

\begin{equation}
... + 9.781133223498595471402477458140584456\text{E}58 x^{45} \;.\label{eq:x^45}
\end{equation}

\section{Conclusions}

We have found the Fourier-Legendre series of modified Bessel functions
of the first kind $I_{N}\left(kx\right)$ based on that found by Keating
\cite{KeatingPhD} for the Bessel functions of the first kind $J_{N}\left(kx\right)$,
and show that Keating's coefficients, comprised of infinite-series,
can be reduced to $\,_{2}F_{3}$ functions. For $N=0$ and $1$ we give
numerical values for those coefficients up through $x^{46}$ with
33-digit accuracy. 

Each of these infinite Fourier-Legendre series may be decomposed into
an infinite sum of infinite series, by gathering like powers from
the Legendre polynomials in each of the terms in the Fourier-Legendre
series. We show that each of these infinite sub-series converges to
values that are inverse powers of the first eight primes $1/\left(2^{i}3^{j}5^{k}7^{l}11^{m}13^{n}17^{o}19^{p}\right)$
multiplying powers of the coefficient $k$. Given the relative paucity of infinite series whose values are known
(e.g., two dozen pages in Gradshteyn and Ryzhik compared to their
900 pages of known integrals), having even one such to add to the
total has the potential to be of use to future researchers. We add
an infinite set of infinite series of $\,_{2}F_{3}$ functions whose
values are now known.

\section*{Appendix}

The code for calculating these series is given in Fortran, below,
because it can be used as such or called in C and C++ programs and
one wishes to avoid duplication in this paper.  %
In Fortran 90, quadruple precision is instituted as in the following
calculation of $\pi$:

\begin{verbatim}
! -----+-----+-----+-----+-----+-----+-----+-----+-----+-----+-----
! In gcc46, and later, one compiles the program with
! $ gfortran precisiontest.f90 -o precisiontest90
! and runs it with
! $ ./precisiontest90
!
! The ouptput is
!                       123456789112345678921234567893123
!                              4       8               16
! s_r_k(6,37)    =    3.1415927    
! s_r_k(15,307)  =    3.1415926535897931     
! s_r_k(33,4931) =    3.1415926535897932384626433832795028 
!
! -----+-----+-----+-----+-----+-----+-----+-----+-----+-----+-----
module precisionkinds!   use ISO_FORTRAN_ENV
   implicit none
   private
   public isp, idp, iqp
   integer, parameter :: isp = selected_real_kind(6, 37)
   integer, parameter :: idp = selected_real_kind(15, 307)
   integer, parameter :: iqp = selected_real_kind(33, 4931)
end module precisionkinds
!
program precisiontest
   use precisionkinds
   implicit none
   real (isp), parameter :: pi1 = 4*atan (1.0_isp)
   real (idp), parameter :: pi2 = 4*atan (1.0_idp)
   real (iqp), parameter :: pi4 = 4*atan (1.0_iqp)
      write (*,*) '                      123456789112345678921234567893123'
      write (*,*) '                             4',                    &
      '       8               16'
          write (*,*) 's_r_k(6,37)    = ',pi1
          write (*,*) 's_r_k(15,307)  = ',pi2
          write (*,*) 's_r_k(33,4931) = ',pi4
end program precisiontest
\end{verbatim}

For readers who are modifying legacy programs, the Fortran 77 equivalent
lines are (two additional compilers are referenced: xlf for PowerPC
and ifort for Intel Macintosh computers):

\begin{verbatim}
C $ xlf -o precisiontest77g precisiontest77g.f
C $ ./precisiontest77g
C $ ifort precisiontest77g.f -o precisiontest77g_intel
C
C                          4       8                 16
C intel gives
Cpi16srk31    =   3.14159265358979323846264338327950      
C
C gfortran gives
C             =   3.14159265358979323846264338327950280 
C xlf gives  
C pi16Q0      =   3.1415926535897932384626433832795059
C pi16srk20   =   3.1415926535897932384626433832795059
C pi16srk31   =   3.1415926535897932384626433832795059
C pi16iesrk31 =   3.1415926535897932384626433832795059
C
C"from Bailey's DQFUN:A thread-safe double-quad precision package
C                 3.141592653589793238462643383279502884197169...
C
          real*16 pi16Q0,pi16srk20,pi16iesrk31,pi16srk31
C
          pi16Q0 = 4*atan (1.0Q0)
          pi16srk20 = 4*atan (Real(1.0,SELECTED_REAL_KIND(20,140)))
          pi16srk31 = 4*atan (Real(1.0,SELECTED_REAL_KIND(31)))
          pi16iesrk31 = 4*atan (Real(1.0,IEEE_SELECTED_REAL_KIND(31)))
\end{verbatim}
In C and C++, one specifies \_Float128. Mathematica and Maple input
may be generated from the following by replacing the Fortran power
operator ``{*}{*}'' by ``\textasciicircum{}'' and deleting the
continuation indicators ``-'' in column 6.

For (\ref{eq:Fourier-Legendre series-numerical}) the code is

\begin{verbatim}
        J0(x) = 0.91973041008976023931442119408062*P(0,x) - 
     -  0.1579420586258518875737139671443637*P(2,x) + 
     -  0.003438400944601109232996887872072915*P(4,x) - 
     -  0.00002919721848828729693660590986125663*P(6,x) + 
     -  1.31735695244778097765561656314328e-7*P(8,x) - 
     -  3.684500844208203027173771096058866e-10*P(10,x) + 
     -  7.011830032993845928208803328211447e-13*P(12,x) - 
     -  9.665964369858912263671995372753346e-16*P(14,x) + 
     -  1.009636276824546446525342170924936e-18*P(16,x) - 
     -  8.266656955927637858991972584174117e-22*P(18,x) + 
     -  5.448244867762758725890082837839430e-25*P(20,x) - 
     -  2.952527182137354751675774606663400e-28*P(22,x) + 
     -  1.338856158858534469080898670096200e-31*P(24,x) - 
     -  5.154913186088512926193234837816582e-35*P(26,x) + 
     -  1.706231577038503450138564028467634e-38*P(28,x) - 
     -  4.906893556427796857473097979568289e-42*P(30,x) + 
     -  1.237489200717479383020539576221293e-45*P(32,x) - 
     -  2.759056237537871868604555688548364e-49*P(34,x) + 
     -  5.477382207172712629199714648396409e-53*P(36,x) - 
     -  9.744200345578852550688946057050674e-57*P(38,x) + 
     -  1.562280711659504489828025148995770e-60*P(40,x) - 
     -  2.269056283827394368836057470594599e-64*P(42,x)
\end{verbatim}
and the corresponding expressions for both forms of (\ref{eq:poly_approx_48_digit})
are

\begin{verbatim}
        J0(x) = 1. - 0.25*x**2 + 0.015625*x**4 - 
     -  0.0004340277777777777777777777777777777778*x**6 + 
     -  6.781684027777777777777777777777777778e-6*x**8 - 
     -  6.781684027777777777777777777777777778e-8*x**10 + 
     -  4.709502797067901234567901234567901235e-10*x**12 - 
     -  2.402807549524439405391786344167296548e-12*x**14 + 
     -  9.385966990329841427311665406903502142e-15*x**16 - 
     -  2.896903392077111551639402903365278439e-17*x**18 + 
     -  7.242258480192778879098507258413196097e-20*x**20 - 
     -  1.496334396734045222954237036862230599e-22*x**22 + 
     -  2.597802772107717400962217077885817011e-25*x**24 - 
     -  3.842903509035084912666001594505646466e-28*x**26 + 
     -  4.901662639075363409012757135849038860e-31*x**28 - 
     -  5.446291821194848232236396817610043178e-34*x**30 + 
     -  5.318644356635593976793356267197307791e-37*x**32 - 
     -  4.600903422695150498956190542558224733e-40*x**34 + 
     -  3.550079801462307483762492702591222788e-43*x**36 - 
     -  2.458504017633176927813360597362342651e-46*x**38 + 
     -  1.5365650110207355798833503733514641567e-49*x**40 - 
     -  8.7106860035189091830121903251216788929e-53*x**42
\end{verbatim}
and (for integer arithmetic),

\begin{verbatim}
J0(x)=x^0*-x^2/(2^2)+x^4/(2^6)-x^6/(2^8*3^2)+x^8/(2^14*3^2)
-x^10/(2^16*3^2*5^2)+x^12/(2^20*3^4*5^2)-x^14/(2^22*3^4*5^2*7^2)
+x^16/(2^30*3^4*5^2*7^2)-x^18/(2^32*3^8*5^2*7^2)+x^20/(2^36*3^8*5^4*7^2)
-x^22/(2^38*3^8*5^4*7^2*11^2)+x^24/(2^44*3^10*5^4*7^2*11^2)
-x^26/(2^46*3^10*5^4*7^2*11^2*13^2)+x^28/(2^50*3^10*5^4*7^4*11^2*13^2)
-x^30/(2^52*3^12*5^6*7^4*11^2*13^2)+x^32/(2^62*3^12*5^6*7^4*11^2*13^2)
-x^34/(2^64*3^12*5^6*7^4*11^2*13^2*17^2)
+x^36/(2^68*3^16*5^6*7^4*11^2*13^2*17^2)
-x^38/(2^70*3^16*5^6*7^4*11^2*13^2*17^2*19^2)
+x^40/(2^76*3^16*5^8*7^4*11^2*13^2*17^2*19^2)
-x^42/(2^78*3^18*5^8*7^6*11^2*13^2*17^2*19^2)
\end{verbatim}

For (\ref{eq:Fourier-Legendre series-numerical_J1})

\begin{verbatim}
        J1(x) = 0.4635981705953810635941110039338702*P(1,x) - 
     -  0.02386534565840739796307209416484866*P(3,x) + 
     -  0.0003197243559720047638524757623256028*P(5,x) - 
     -  1.970519180666594250258062929391112e-6*P(7,x) + 
     -  6.987247473097807218791759410157014e-9*P(9,x) - 
     -  1.610500056046875027807002442953327e-11*P(11,x) + 
     -  2.607086592441628842939248193619909e-14*P(13,x) - 
     -  3.127311482540796882144713619567442e-17*P(15,x) + 
     -  2.891424081787050739827382596616064e-20*P(17,x) - 
     -  2.123664534779369199214414455720317e-23*P(19,x) + 
     -  1.269011201758673511714553707528186e-26*P(21,x) - 
     -  6.290201939135925763576871358738600e-30*P(23,x) + 
     -  2.628135796989325452573870774267213e-33*P(25,x) - 
     -  9.381575562723076109283258050667642e-37*P(27,x) + 
     -  2.894337242415984040941859061022419e-40*P(29,x) - 
     -  7.794444104104171684395094261174814e-44*P(31,x) + 
     -  1.848200759818170134895867052306767e-47*P(33,x) - 
     -  3.888249639773912225694535890329244e-51*P(35,x) + 
     -  7.306978718807123633044120058516188e-55*P(37,x) - 
     -  1.234022530456621571127590099647796e-58*P(39,x) + 
     -  1.883067799255568915649461884255428e-62*P(41,x) - 
     -  2.609122884536350861268195351045890e-66*P(43,x)
\end{verbatim}
and the corresponding expressions for both forms of (\ref{eq:J1series})
are

\begin{verbatim}
        J1(x)= 0.5*x - 0.0625*x**3 + 
     -  0.002604166666666666666666666666666666667*x**5 - 
     -  0.00005425347222222222222222222222222222222*x**7 + 
     -  6.781684027777777777777777777777777778e-7*x**9 - 
     -  5.651403356481481481481481481481481481e-9*x**11 + 
     -  3.363930569334215167548500881834215168e-11*x**13 - 
     -  1.501754718452774628369866465104560343e-13*x**15 + 
     -  5.214426105738800792950925226057501190e-16*x**17 - 
     -  1.448451696038555775819701451682639219e-18*x**19 + 
     -  3.291935672814899490499321481096907317e-21*x**21 - 
     -  6.234726653058521762309320986925960827e-24*x**23 + 
     -  9.991549123491220772931604145714680813e-27*x**25 - 
     -  1.372465538941101754523571998037730881e-29*x**27 + 
     -  1.633887546358454469670919045283012953e-32*x**29 - 
     -  1.701966194123390072573874005503138493e-35*x**31 + 
     -  1.564307163716351169645104784469796409e-38*x**33 - 
     -  1.278028728526430694154497372932840204e-41*x**35 + 
     -  9.342315267006072325690770269976902073e-45*x**37 - 
     -  6.146260044082942319533401493405856627e-48*x**39 + 
     -  3.658488121477941856865119936551105135e-51*x**41 - 
     -  1.979701364436115723411861437527654294e-54*x**43
\end{verbatim}

and (for integer arithmetic),

\begin{verbatim}
J1(x)=x^2/2-x^3/(2^4)+x^5/(2^7*3*-x^7/(2^11*3^2)+x^9/(2^15*3^2*5)
-x^11/(2^18*3^3*5^2)+x^13/(2^21*3^4*5^2*7)-x^15/(2^26*3^4*5^2*7^2)
+x^17/(2^31*3^6*5^2*7^2)-x^19/(2^34*3^8*5^3*7^2)+x^21/(2^37*3^8*5^4*7^2*11)
-x^23/(2^41*3^9*5^4*7^2*11^2)+x^25/(2^45*3^10*5^4*7^2*11^2*13)
-x^27/(2^48*3^10*5^4*7^3*11^2*13^2)+x^29/(2^51*3^11*5^5*7^4*11^2*13^2)
-x^31/(2^57*3^12*5^6*7^4*11^2*13^2)+x^33/(2^63*3^12*5^6*7^4*11^2*13^2*17)
-x^35/(2^66*3^14*5^6*7^4*11^2*13^2*17^2)
+x^37/(2^69*3^16*5^6*7^4*11^2*13^2*17^2*19)
-x^39/(2^73*3^16*5^7*7^4*11^2*13^2*17^2*19^2)
+x^41/(2^77*3^17*5^8*7^5*11^2*13^2*17^2*19^2)
-x^43/(2^80*3^18*5^8*7^6*11^3*13^2*17^2*19^2)
\end{verbatim}

For (\ref{eq:Fourier-Legendre seriesI0}) we have

\begin{verbatim}
         I(0,x) = 1.086521097023589815837941923492506*P(0,x) + 
     -  0.1758046819215242662605951354261250*P(2,x) + 
     -  0.003709009244052882533923838165527033*P(4,x) + 
     -  0.00003095105270992432198613744608777602*P(6,x) + 
     -  1.381259734719773538320052305224506e-7*P(8,x) + 
     -  3.834312601086373005317788906125573e-10*P(10,x) + 
     -  7.257172450096213936720667660411978e-13*P(12,x) + 
     -  9.962746978836018020128433111635975e-16*P(14,x) + 
     -  1.037251346110052630963705477046736e-18*P(16,x) + 
     -  8.470496863240475339343499321604116e-22*P(18,x) + 
     -  5.5705413998588522192782604835236869e-25*P(20,x) + 
     -  3.013347383234528850224689041823201e-28*P(22,x) + 
     -  1.364338005353527272638479093175249e-31*P(24,x) + 
     -  5.246088467162281648944660989359565e-35*P(26,x) + 
     -  1.734417052236546525979562610169336e-38*P(28,x) + 
     -  4.982929889631203560686967762401821e-42*P(30,x) + 
     -  1.255546430559877621587201790700357e-45*P(32,x) + 
     -  2.797096596401706413444068821508193e-49*P(34,x) + 
     -  5.548955677049963483909673845071489e-53*P(36,x) + 
     -  9.865206225205083247212985096531573e-57*P(38,x) + 
     -  1.580763691652306983099443761944673e-60*P(40,x) + 
     -  2.294688331479205281600814719914093e-64*P(42,x) + 
     -  3.031771495580703895127109933386607e-68*P(44,x) + 
     -  3.661200772680598752990852186025167e-72*P(46,x)
\end{verbatim}

and for (\ref{eq:Fourier-Legendre series-I1}) we have

\begin{verbatim}
        I(1,x) = 0.5386343421852555592809081051666336*P(1,x) + 
     -  0.02618069164825977449795296407260333*P(3,x) + 
     -  0.0003419851912550806236210094361507344*P(5,x) + 
     -  2.077651971699656963860267070724864e-6*P(7,x) + 
     -  7.299001518662431414905576324932877e-9*P(9,x) + 
     -  1.671443482954853739162527767203215e-11*P(11,x) + 
     -  2.692744551459235232734936666452704e-14*P(13,x) + 
     -  3.218106754771162455853759838545282e-17*P(15,x) + 
     -  2.966624646773403074824196435937542e-20*P(17,x) + 
     -  2.173686883720901031568436047655748e-23*P(19,x) + 
     -  1.296326265789554875546711294998344e-26*P(21,x) + 
     -  6.414855102151415733596588296578833e-30*P(23,x) + 
     -  2.676389925142875786863074285387196e-33*P(25,x) + 
     -  9.542035089444710700263714658817980e-37*P(27,x) + 
     -  2.940669572337884276201779460203862e-40*P(29,x) + 
     -  7.911705033029330504663434638574930e-44*P(31,x) + 
     -  1.874426565726980007813411585706840e-47*P(33,x) + 
     -  3.940459072597980181771454632976205e-51*P(35,x) + 
     -  7.400090413796917559009186360838868e-55*P(37,x) + 
     -  1.248984620737396858084740490332061e-58*P(39,x) + 
     -  1.904842982553207494042180613785837e-62*P(41,x) + 
     -  2.637959760920312924684635466402215e-66*P(43,x) + 
     -  3.332061910821697596383220274010501e-70*P(45,x)
\end{verbatim}

\end{document}